\documentclass[letterpaper, 10 pt, conference]{ieeeconf}

\IEEEoverridecommandlockouts
\usepackage{graphics}
\usepackage{epsfig}
\usepackage{mathptmx}
\usepackage{algorithm}
\usepackage{algorithmic}
\usepackage{times}
\usepackage{amsmath}
\usepackage{amssymb}
\usepackage{dsfont}
\usepackage{cleveref}
\usepackage{booktabs}
\usepackage{tabularx}
\usepackage{ragged2e}
\usepackage{mathtools}
\newlength{\opws}
\newlength{\opwl}
\newcommand{\opS}[1]{\mathmakebox[\opws][l]{#1}}
\newcommand{\opL}[1]{\mathmakebox[\opwl][l]{#1}}
\newcommand{\opM}[1]{\mathop{\mathmakebox[\opwl][l]{\min}}\limits_{\mathmakebox[\opwl][l]{#1}}}
\usepackage{tikz}
\usetikzlibrary{arrows.meta}
\usepackage{cite}

\newcommand{\PP}{\mathbb{P}}
\newcommand{\ind}[1]{\mathds{1}_{#1}}
\newcommand{\Jhat}{\widehat{J}}
\newcommand{\Jec}{J_{\text{ec}}}
\newcommand{\Nmean}{N}
\newcommand{\Nrisk}{N_{\text{SAA}}^{\text{risk}}}
\newcommand{\CVaR}{\text{CVaR}}

\usepackage{xcolor}
\usepackage{url}

\title{\LARGE \bf
Risk Control of Traffic Flow Through Chance Constraints and Large Deviation Approximation
}

\author{Rui Xu$^{1}$, Shanyin Tong$^{2}$ and Xuan Di$^{1,3}$%
\thanks{This work was supported by NSF under Grant DMS-2529292.}%
\thanks{$^1$Department of Civil Engineering and Engineering Mechanics,
        Columbia University, New York,
NY 10027 USA
        {\tt\small (Emails: rx2291@columbia.edu; sharon.di@columbia.edu)}}%
\thanks{$^2$Department of Mathematics, University of Pennsylvania, Philadelphia, PA 19104 USA
        {\tt\small (Email: tong3@sas.upenn.edu)}}%
\thanks{$^3$Data Science Institute, Columbia University, New York, NY 10027 USA}%
}

\begin{document}

\maketitle
\thispagestyle{empty}
\pagestyle{empty}

\begin{abstract}

We develop a rare chance-constrained optimal control framework for autonomous traffic management, with traffic dynamics governed by stochastic partial differential equations.
To efficiently enforce these probabilistic safety specifications, we exploit a large deviation theory (LDT) based approximation method, which converts the original highly non-convex, sampling-heavy stochastic optimization problem into a tractable deterministic nonlinear programming problem. 
The proposed LDT-based reformulation exhibits superior computational scalability, as its computational cost remains approximately constant regardless of the target violation probability, as opposed to the scaling challenge of traditional sampling-based methods. 
The effectiveness of the proposed framework %
is illustrated through numerical experiments on 
a ring road with stochastic on-ramp inflow, 
and evaluated with two representative risk measures, namely, maximum traffic density and speed variation. 
The results show that the proposed method keeps the violation probability close to the target without parameter tuning and is more than ten times faster than a risk-averse baseline.

\end{abstract}

\section{INTRODUCTION}

Uncertainties in traffic systems, such as fluctuating on-ramp demand and upstream inflow, degrade traffic efficiency and safety. They produce excessive density and large speed variations,  which are closely related to crash risks \cite{lee2002crash}. 
Despite these significant impacts, commonly used risk-neutral control methods aim to optimize expected values under uncertainty \cite{chu2019cacc}, %
thus focusing on average traffic behavior without explicitly accounting for these risks.
This motivates us to develop risk-aware traffic control methods accounting for the  tail behavior (extreme outcomes) under uncertainty.

Risk-aware control methods that account for traffic risk can be divided into two bodies of literature: objective-based approaches that penalize risk through the objective function, and constraint-based approaches that explicitly impose risk as additional constraints. 
Objective-based methods are implemented either by evaluating the objective under the worst-case scenario \cite{shang2024robust}, or by directly penalizing tail behavior in the objective, as in conditional value-at-risk (CVaR) control, which penalizes the average cost over the worst fraction of uncertainty realizations \cite{shi2024perimeter}.
Unlike objective-based methods, which account for risk indirectly through penalties, 
constraint-based methods impose explicit safety requirements and provide direct control over undesirable events, making them more suitable for traffic applications %
\cite{zhao2020drcc_cacc}.
These constraints can be imposed deterministically with control barrier functions (CBFs) that maintain the system within a safe set for each prescribed disturbance realization \cite{ames2019cbf,alan2023cbf_av}. 
Alternatively, they can be imposed probabilistically by requiring the probability of an extreme event to stay below a prescribed risk threshold, leading to standard chance constraints (CCs) under a given uncertainty distribution \cite{wu2010chance_ramp,lyu2021merging}, and distributionally robust chance constraints (DRCCs) over a family of plausible distributions \cite{zhao2020drcc_cacc}. Because chance constraints neither rely on a predefined set of realizations nor require adding a large number of realization-specific constraints, they are well suited to traffic problems and are the focus of this work.
Despite these advantages, CCs and DRCCs remain computationally challenging. Sampling-based reformulations of CCs require increasingly many samples as the target probability decreases, with each sample requiring a solution of the traffic model \cite{mehr2017stl,tong2022rare}. DRCCs are generally difficult to solve, and existing tractable reformulations rely on special assumptions, such as convexity or affine dependence on the uncertainty \cite{zhao2020drcc_cacc}. These assumptions are not satisfied by the nonlinear traffic risk maps considered here. \Cref{tab:risk_compare} summarizes representative risk-aware traffic control studies and their treatments of uncertainty, risk, and control.

\begin{table*}[htbp]
\vspace*{5pt}
\centering
\caption{Representative risk-aware traffic control methods}
\label{tab:risk_compare}
\footnotesize
\setlength{\tabcolsep}{3pt}
\renewcommand{\arraystretch}{1.1}

\begin{tabularx}{\textwidth}{@{} 
  >{\RaggedRight\arraybackslash}p{0.12\textwidth} 
  >{\RaggedRight\arraybackslash}p{0.09\textwidth} 
  >{\RaggedRight\arraybackslash}X 
  >{\RaggedRight\arraybackslash}X 
  >{\RaggedRight\arraybackslash}X 
  >{\RaggedRight\arraybackslash}X 
@{}}
\toprule
\textbf{Reference} & \textbf{Method} & \textbf{Traffic model} & \textbf{Risk measure} & \textbf{Decision variable} & \textbf{Methodology} \\
\midrule

Zhao \& Zhang~\cite{zhao2020drcc_cacc}
& DRCC
& Car-following vehicle kinematics
& Worst-case probability of unsafe spacing
& Connected and automated vehicle acceleration
& Semidefinite program via strong duality \\ 
\addlinespace[3pt]

Wu et al.~\cite{wu2010chance_ramp}
& CC
& Zone-based flow conservation
& Probability of flow exceeding capacity
& Ramp metering rate
& Deterministic constraint via Gaussian quantile \\ 
\addlinespace[3pt]

Shi et al.~\cite{shi2024perimeter}
& CVaR
& Macroscopic fundamental diagram
& CVaR of total time spent
& Perimeter transfer-flow proportion
& Convex scenario-tree program via duality \\ 
\addlinespace[3pt]

Shang et al.~\cite{shang2024robust}
& Min--max robust
& Data-driven mixed-traffic car-following model
& Worst-case cost of velocity and spacing errors
& Connected and automated vehicle acceleration
& Quadratic program via vertex or duality methods \\ 
\addlinespace[3pt]

Lyu et al.~\cite{lyu2021merging}
& CC-CBF
& Point-mass vehicle kinematics
& Probability of violating the safe-distance CBF
& Autonomous vehicle acceleration
& Quadratic program via Gaussian quantile \\

\midrule

\textbf{This paper}
& \textbf{Rare CC}
& \textbf{Stochastic PDE}
& \textbf{Probability of extreme density or speed variation}
& \textbf{Autonomous vehicle velocity field}
& \textbf{Nonlinear program via sampling-free LDT} \\

\bottomrule
\end{tabularx}
\vspace{-5pt}
\end{table*}

Motivated by these gaps, this paper develops a rare chance-constrained optimal control framework for autonomous traffic management. The traffic dynamics follow a stochastic macroscopic partial differential equation (PDE), and the control is the velocity field of the autonomous vehicles (AVs). We estimate the probability with large deviation theory (LDT), which represents the rare event probability through the most probable point (MPP), the most likely realization of the extreme traffic event, and computes the MPP from an optimization problem \cite{tong2021extreme,tong2022rare}. Embedding this estimate in the control problem turns the stochastic PDE-constrained problem into a deterministic nonlinear program. The size of this program is independent of the target probability, and the program can be efficiently solved by gradient-based methods.

The contributions of this paper include: 
(1) The first integration of a traffic flow control framework governed by stochastic PDEs with rare chance constraints to regulate the probabilities of extreme traffic events; 
(2) A sampling-free reformulation of the rare chance constraint through LDT, which turns the stochastic PDE-constrained problem into a deterministic nonlinear program whose computational complexity does not scale with the rareness of the event; 
and (3) Extensive numerical validation demonstrating that the proposed method achieves precise near-target probability control without heuristic parameter tuning, outperforming conventional risk-averse baselines in both reliability and computational efficiency.

The remainder of this paper is organized as follows. Section~II presents the stochastic traffic flow control model, the traffic risk measures, and the LDT-based reformulation. Section~III presents the numerical experiments and compares the proposed method with three baselines. Section~IV concludes the paper and discusses future directions.

\section{Methodology}

\subsection{PDE-constrained stochastic optimal control}

We will first present a macroscopic traffic flow model to avoid controlling each vehicle individually, 
which is supported by literature on mean field games \cite{huang2020mfg} and on traffic PDE models \cite{yu2021output}.

The objective of the stochastic traffic flow control problem is to find an optimal autonomous vehicle velocity field $u(x,t)$ that minimizes the expected cost, subject to the traffic flow dynamics following a stochastic forced continuity equation (CE):
\begin{equation}
\label{eq:field_control}
\begin{aligned}
    \opS{\min\limits_{u}}\quad & \mathbb{E}\big[J(u,\rho)\big] \\
    \opS{\text{s.t.}}\quad
    &  \rho_t + (\rho u)_x = s(x,t;\theta), \quad \theta \sim \pi,
\end{aligned}
\end{equation}
where $J(u,\rho)$ denotes the total cost depending on the traffic density $\rho(x,t)$ and the velocity $u$, and $s(x,t;\theta)$ is the stochastic source term modeling demand perturbations with $\theta$ representing the uncertainty parameter and following a given distribution $\pi$. The control $u(x,t)$ is the actual speed of the AVs. The system also satisfies the initial condition $\rho(x,0) = \rho_0(x)$ and the periodic boundary condition $\rho(0,t) = \rho(L,t)$ over the space-time domain $(x,t) \in [0,L] \times [0,T]$, which will be omitted in the following sections for simplicity. %

In reality, the uncertainty is often introduced by on-ramp flows or demand fluctuations. Specifically, we use a masked log-normal random field based on the Karhunen--Lo\`eve (KL) expansion \cite{betz2014numerical} to model the stochastic source term: 
\begin{equation}
\label{eq:kl}
    s(x,t;\theta)
    =
        \left(\mu_s
        +
        \sum_{k=1}^{d}
        \sqrt{\lambda_k}\, f_k(x)\exp(\theta_k) \right)\cdot\ind{\mathcal{D}_{\mathrm{mask}}}(x,t),
\end{equation}
where $\mathcal{D}_{\mathrm{mask}} \subset [0, L]\times[0, T]$ represents the localized space-time domain of the on-ramp. The scalar $\mu_s$ is the mean source strength, $(\lambda_k, f_k(x))$ are eigenpairs, obtained analytically, for the covariance kernel $C(x,y) = \sigma^2 \exp(-|x-y|/l_c)$ with zero boundary conditions at both endpoints, and the variance $\sigma^2$ and the correlation length $l_c$ are fixed. The underlying random vector $\theta\in\mathbb{R}^d$ follows a multivariate Gaussian distribution $\theta \sim \pi = \mathcal{N}(\mu_\theta,\Sigma_\theta)$ with mean $\mu_\theta$ and covariance $\Sigma_\theta$.

\subsection{Traffic risk measure}\label{sec:risk-measure}

While the expected objective in \eqref{eq:field_control} regulates average performance, our interest here is to explicitly constrain the occurrence of extreme events (such as collisions). 
Thus, we first discuss the risk measure. %

We introduce a traffic risk measure operator $M$ and define the extreme-event map $F$ through its maximum. 
Given a selected measure operator $M(\rho, u)$, we define the corresponding extreme-event map by
\begin{equation}\label{eq:forward_max}
    F(u,\theta)
    :=\max_{(x,t)\in[0,L]\times[0,T]} M(\rho(x,t), u(x,t)),
\end{equation}
where the dependence on $u$ and $\theta$ is induced through the state $\rho$ solving the CE in \eqref{eq:field_control}.

In this work, we consider two representative choices of risk measure $M(\rho(x,t), u(x,t))$:

\subsubsection{Maximum density}
$M(\rho(x,t), u(x,t))=\rho(x,t)$. Bounding the maximum density serves as a direct macroscopic surrogate for collision avoidance, as it mathematically prevents the average inter-vehicle spacing from falling below the physical vehicle length ($1/\rho(x,t) \le L_{\mathrm{veh}}$).

\subsubsection{Maximum coefficient of variation of speed (CVS)} $M(\rho(x,t), u(x,t)) = \sqrt{\mathrm{Var}(u)(x,t)} / \bar{u}_\Delta(x,t)$, capturing the local speed variation (i.e., the relative standard deviation of the speed $u$) over a spatial window spanning $\Delta$ grids. Here, the local mean speed is defined as $\bar{u}_\Delta(x,t) = M_{\Delta}(x,t)^{-1} \int_{x}^{x+\Delta} \rho(y,t)u(y,t) \, dy$ and the local speed variance is $\mathrm{Var}(u)(x,t) = M_{\Delta}(x,t)^{-1} \int_{x}^{x+\Delta} \rho(y,t)(u(y,t) - \bar{u}_\Delta(x,t))^2 \, dy$, with the local mass $M_{\Delta}(x,t) = \int_{x}^{x+\Delta} \rho(y,t) \, dy$. This measure captures instability that may arise from strong local speed variation even before density itself becomes critical.

\subsection{Traffic risk control through rare chance constraints}\label{sec:cc}
To explicitly enforce that the probability of traffic risk remains low (below a given risk threshold), we add a chance constraint to the stochastic control problem \eqref{eq:field_control} and formulate the following rare chance-constrained stochastic control problem:
\begin{subequations}\label{eq:cc_problem}
\begin{align}
    \opS{\min\limits_{u}}\quad & \mathbb{E}\big[J(u,\rho)\big] \label{eq:cc_problem_obj}\\
    \opS{\text{s.t.}}\quad
    &  \rho_t + (\rho u)_x = s(x,t;\theta), \quad \theta \sim \pi, \label{eq:cc_problem_pde}\\
    & \mathbb{P}\big(F(u,\theta)\ge z\big)\le \alpha, \label{eq:cc_problem_cc}
\end{align}
\end{subequations}
where $z$ is the prescribed safety threshold and $\alpha\ll 1$ is the allowable violation probability.

The rare chance constraint \eqref{eq:cc_problem_cc} is highly non-convex and computationally expensive to evaluate. Sampling-based reformulations, such as the scenario approach and sample average approximation, require a number of samples that grows in proportion to $1/\alpha$, and each sample enters the optimization problem through one extra constraint requiring one solution of the traffic model, which makes the problem large-scale in the rare-event regime \cite{tong2022rare}. To overcome this, we utilize a sampling-free approximation scheme based on LDT for rare chance-constrained optimization \cite{tong2021extreme,tong2022rare}.

\subsection{LDT reformulation of chance-constrained optimization}

Before explaining how the LDT reformulation works, we first write \eqref{eq:cc_problem} in a more abstract form. 
For a given control $u$ and parameter $\theta$, the continuity equation in \eqref{eq:field_control} has the solution $\rho(u,\theta)$. Denote $\Jhat(u,\theta):=J(u,\rho(u,\theta))$ for the cost at that solution, and keep the same way as we define the extreme-event map $F(u,\theta)$ in \eqref{eq:forward_max}. The objective and the chance constraint then depend on the control and the uncertain parameter alone,
\begin{subequations}\label{eq:cc_reduced}
\begin{align}
    \opS{\min\limits_{u}}\quad & \mathbb{E}\big[\Jhat(u,\theta)\big] \label{eq:cc_reduced_obj}\\
    \opS{\text{s.t.}}\quad & \mathbb{P}\big(F(u,\theta)\ge z\big)\le \alpha. \label{eq:cc_reduced_cc}
\end{align}
\end{subequations}

Under certain regularity assumptions on $F$, LDT connects the rare event probability to an optimization problem,
\begin{equation}\label{eq:ldt}
  \mathbb{P}\big(F(u,\theta)\ge z\big)
  \asymp
  \exp\!\Big(-\tfrac{1}{2}\big\|\theta^\ast(u,z)-\mu_\theta\big\|_{\Sigma_\theta^{-1}}^{2}\Big),\;\; \text{as } z\to\infty,
\end{equation}
where $\asymp$ denotes the log-asymptotics\footnote{$\asymp$ means the ratio between the logarithms of LHS and RHS goes to one as $z \to\infty$.} and $\theta^\ast$ is the LDT minimizer, also called the instanton; in the Gaussian case, it is the most probable point (MPP). The LDT minimizer $\theta^\ast$ minimizes the rate function over the rare event region and is obtained by solving the following optimization problem
\begin{equation}
\label{eq:mpp}
    \theta^*(u, z):=\underset{\theta \in \mathbb{R}^d}{\arg\,\min} \; \left\{ \frac{1}{2}\|\theta-\mu_\theta\|_{{\Sigma}_\theta^{-1}}^2 : F(u,\theta)\ge z \right\},
\end{equation}
where $\|x\|_{\Sigma_\theta^{-1}}:=(x^\top\Sigma_\theta^{-1}x)^{1/2}$ and the rate function of the Gaussian $\theta$ is the squared $\Sigma_\theta^{-1}$ norm.

Following the LDT framework \cite{tong2021extreme,tong2022rare,schorlepp2023scalable}, the rare-event probability is approximated by the measure of the half-space determined by the LDT minimizer (also called the first-order reliability method (FORM) \cite{du2001most}):

\begin{equation}\label{eq:form_prob}
    \mathbb{P}\big(F(u,\theta)\ge z\big) \approx \Phi\left( -\|\theta^\ast-\mu_\theta\|_{{\Sigma}_\theta^{-1}} \right),
\end{equation}
where $\Phi(\cdot)$ is the cumulative distribution function of the standard normal distribution.

Substituting \eqref{eq:form_prob} into the chance constraint \eqref{eq:cc_reduced_cc} gives a bilevel problem. The upper level optimizes the control $u$, and the lower level \eqref{eq:mpp} returns the MPP $\theta^\ast$ that enters the upper-level constraint.

\noindent\textbf{[Bilevel problem]}
\begin{subequations}\label{eq:bilevel}
\begin{alignat}{3}
  & {\min_{u}} \quad && \mathbb{E}\big[\Jhat(u,\theta)\big] && \label{eq:bilevel_obj}\\
  & \text{s.t.}\quad && \Phi\big(-\|\theta^\ast-\mu_\theta\|_{\Sigma_\theta^{-1}}\big)\le \alpha, && \label{eq:bilevel_prob}\\
  &&& \theta^\ast=\underset{\theta \in \mathbb{R}^d}{\arg\,\min} \; \left\{ \frac{1}{2}\|\theta-\mu_\theta\|_{{\Sigma}_\theta^{-1}}^2 : F(u,\theta)\ge z \right\}. && \label{eq:bilevel_lower}
\end{alignat}
\end{subequations}

\begin{figure*}[tb]
\vspace*{5pt}
\centering
\resizebox{\textwidth}{!}{%
\begin{tikzpicture}[
  box/.style={draw=black!60, rounded corners=2.5pt, very thick, inner sep=6pt},
  ar/.style={-{Latex[length=2.6mm]}, very thick}]
\node[box] (b) at (0,0) {\begin{minipage}[c][2.6cm][t]{6.2cm}\centering
  \textbf{Reduced problem}~\eqref{eq:cc_reduced}\\[0.6em]
  $\begin{array}{@{}l@{\;}l@{}}
     \opL{\min\limits_{u}} & \mathbb{E}\big[\Jhat(u,\theta)\big]\\[0.35em]
     \opL{\text{s.t.}} & \mathbb{P}\big(F(u,\theta)\ge z\big)\le\alpha
   \end{array}$
  \end{minipage}};
\node[box] (c) at (7.95,0) {\begin{minipage}[c][2.6cm][t]{6.2cm}\centering
  \textbf{Bilevel problem}~\eqref{eq:bilevel}\\[0.6em]
  $\begin{array}{@{}l@{\;}l@{}}
     \opL{\min\limits_{u}} & \mathbb{E}\big[\Jhat(u,\theta)\big]\\[0.35em]
     \opL{\text{s.t.}} & \Phi\big(-\|\theta^\ast-\mu_\theta\|_{\Sigma_\theta^{-1}}\big)\le\alpha\\[0.35em]
     & \theta^\ast=\begin{array}[t]{@{}l@{\;}l@{}}
                     \arg\min\limits_{\theta} & \tfrac12\|\theta-\mu_\theta\|^2_{\Sigma_\theta^{-1}}\\[0.25em]
                     \text{s.t.} & F(u,\theta)\ge z
                   \end{array}
   \end{array}$
  \end{minipage}};
\node[box] (d) at (15.9,0) {\begin{minipage}[c][2.6cm][t]{6.2cm}\centering
  \textbf{Single-level NLP}~\eqref{eq:nlp}\\[0.6em]
  $\begin{array}{@{}l@{\;}l@{}}
     \opM{u,\theta^\ast,\lambda\ge0} & \tfrac{1}{\Nmean}\sum\nolimits_{i=1}^{\Nmean} \Jhat(u,\theta^{(i)})\\[0.55em]
     \opL{\text{s.t.}} & \|\theta^\ast-\mu_\theta\|_{\Sigma_\theta^{-1}}\ge\Phi^{-1}(1-\alpha)\\[0.35em]
     & F(u,\theta^\ast)=z\\[0.35em]
     & \Sigma_\theta^{-1}(\theta^\ast-\mu_\theta)=\lambda\nabla_\theta F(u,\theta^\ast)
   \end{array}$
  \end{minipage}};
\draw[ar] (b.east) -- node[above,align=center,font=\small\bfseries]{LDT \\ FORM} (c.west);
\draw[ar] (c.east) -- node[above,align=center,font=\small\bfseries]{KKT\\SAA} (d.west);
\end{tikzpicture}}
\caption{Flow of LDT reformulation of the reduced problem \eqref{eq:cc_reduced} into the deterministic nonlinear program \eqref{eq:nlp}. The first arrow replaces the probability by the LDT approximation/FORM \eqref{eq:form_prob}. The second arrow replaces the lower-level problem by its KKT conditions and the expectation in the objective by a sample average.}
\label{fig:pipeline}
\vspace{-5pt}
\end{figure*}

We transform \eqref{eq:bilevel} line by line into a deterministic single-level nonlinear program (NLP).
The objective \eqref{eq:nlp_obj} is the sample average approximation (SAA) \cite{kleywegt2002sample} of \eqref{eq:bilevel_obj} over $\Nmean$ independent and identically distributed (i.i.d.) samples $\theta^{(i)}\sim\pi$, and each sample requires one solution of the traffic flow model.
Constraints \eqref{eq:nlp_limit} and \eqref{eq:nlp_kkt} replace the lower-level problem \eqref{eq:bilevel_lower} by its Karush--Kuhn--Tucker (KKT) conditions.
Constraint \eqref{eq:nlp_limit} places the MPP $\theta^\ast$ on the boundary of the rare event region, and \eqref{eq:nlp_kkt} is the stationarity condition, where $\lambda\ge 0$ is the Lagrange multiplier of the constraint $F(u,\theta)\ge z$ in \eqref{eq:bilevel_lower}.
The resulting NLP is as follows.

\noindent\textbf{[Single-level NLP]}
\begin{subequations}\label{eq:nlp}
\begin{alignat}{3}
  & {\min_{\mathclap{u,\,\theta^\ast,\,\lambda\ge 0}}} \qquad && \frac{1}{\Nmean}\sum_{i=1}^{\Nmean} \Jhat(u,\theta^{(i)}) && \label{eq:nlp_obj}\\
  & \text{s.t.}\qquad && \|\theta^\ast-\mu_\theta\|_{\Sigma_\theta^{-1}} \ge \Phi^{-1}(1-\alpha), && \label{eq:nlp_prob}\\
  &&& F(u,\theta^\ast)=z, && \label{eq:nlp_limit}\\
  &&& \Sigma_\theta^{-1}(\theta^\ast-\mu_\theta) = \lambda \nabla_\theta F(u,\theta^\ast). && \label{eq:nlp_kkt}
\end{alignat}
\end{subequations}

This single-level deterministic NLP formulation eliminates the need for large-scale Monte Carlo sampling to evaluate the rare-event chance constraint by combining the LDT approximation with the KKT relaxation in \eqref{eq:nlp_limit}--\eqref{eq:nlp_kkt}, and it can be solved efficiently using off-the-shelf gradient-based optimization solvers (e.g., IPOPT \cite{wachter2006implementation}).

\Cref{fig:pipeline} summarizes the flow of the LDT reformulation. The first step obtains the reduced problem \eqref{eq:cc_reduced} from \eqref{eq:cc_problem} by incorporating the PDE implicitly into the objective and risk measure.
The second step replaces the probability in the constraint by its LDT approximation \eqref{eq:form_prob} and gives the bilevel problem \eqref{eq:bilevel}.
The third step replaces the lower-level problem by its KKT conditions and the expectation in the objective by a sample average and gives the single-level program \eqref{eq:nlp}.

\section{Numerical Experiments}

In this section, we study the performance of the proposed LDT-based chance-constrained (LDT-CC) method for risk control of traffic flows by comparing it with three baseline frameworks: the unconstrained, nonseparable, and CVaR-penalized optimal control problems. In all numerical experiments, the fundamental objective is the equilibrium cost, defined as
\begin{equation}
\Jec(u,\rho)
:=
\int_0^T\!\!\int_0^L
\frac12 \rho(x,t)\big(U(\rho)-u\big)^2\,dx\,dt,
\label{eq:J_ec_num}
\end{equation}
which penalizes deviations from the Greenshields equilibrium velocity $U(\rho) = 1-\rho$. This form differs from the usual Greenshields relation because the density is normalized by the jam density and the speed by the free-flow speed, so that $\rho, u, U(\rho)\in[0,1]$. The three baselines are formulated as follows: (1) \emph{Unconstrained} optimal control problem: minimizes the expectation of the equilibrium cost $\mathbb{E}[\Jec(u,\rho)]$. (2) \emph{Nonseparable} optimal control problem: minimizes a composite cost that 
induces vehicles to decelerate in high-density areas
and accelerate in low-density regions \cite{huang2020mfg}:
\begin{equation}
\mathbb{E}\!\big[\Jec(u,\rho)\big]
-
\mathbb{E}\!\left[
\int_0^T\!\!\int_0^L
\frac12 \rho \cdot U(\rho)^2\,dx\,dt
\right].
\label{eq:J_nonsep_num}
\end{equation}
(3) The \emph{CVaR}-penalized optimal control problem adds a penalty for the risk measure to the objective. This is achieved by penalizing $\CVaR_{1-\alpha}(F-z)$, i.e., the average exceedance value of $F-z$ in the worst $\alpha$ fraction of outcomes. In the Rockafellar and Uryasev format \cite{kouri2016cvar,rockafellar2002conditional}, this can be written as
\begin{equation}
\min_{u,\tau}\;
\mathbb{E}[\Jec(u,\rho)]
+
\eta\!\left(
\tau+\frac{1}{\alpha}\,
\mathbb{E}\big[(F(u,\theta)- z -\tau)_+\big]
\right),
\label{eq:cvar_num}
\end{equation}
where $\eta>0$ is the penalty weight and $\tau$ is an auxiliary value-at-risk variable. Our proposed \emph{LDT-CC} method solves the nonlinear program with deterministic constraints as shown in \eqref{eq:nlp}.

The continuous macroscopic PDE is discretized using an explicit finite difference method with Lax--Friedrichs flux and forward Euler time integration (CFL = $0.5$). Implemented in Python, the framework uses JAX for exact automatic differentiation (backward gradients) and IPOPT with the MA57 linear solver \cite{wachter2006implementation}. Both CVaR and LDT-CC are initialized using the unconstrained solution. For gradient-based solvability, the extreme-event and positive part operators are smoothly approximated with fixed parameters $\gamma=\epsilon=0.01$: $F(u, \theta) = \gamma \log \big((LT)^{-1} \int_0^T \int_0^L \exp(M(\rho,u)/\gamma) dxdt\big) \approx \max_{(x,t)} M(\rho,u)$ and $(y)_+ \approx \epsilon \log(1 + \exp(y/\epsilon))$.

Across all cases, the stochastic source $s(x,t;\theta)$ is realized via a 10-mode truncated KL expansion of a log-normal random field ($\mu_s=0.1, \sigma=1.0, l_c=0.5$), defined in \eqref{eq:kl}. The underlying Gaussian vector is parameterized by mean ${\mu}_\theta = -0.5 \mathbf{1}$ and covariance matrix ${\Sigma}_\theta = 0.1 \mathbf{I}$, where $\mathbf{1}$ and $\mathbf{I}$ denote the all-ones vector and identity matrix of appropriate dimensions, respectively. 
We use different sample sizes $\Nmean=100$ and $\Nrisk=5000$ ($\Nrisk=12500$ for $\alpha=0.2\%$) for the expectations of the equilibrium cost and the CVaR risk term, respectively. Since CVaR averages over the worst $\alpha$ fraction of the samples, we scale $\Nrisk$ with $1/\alpha$ so that approximately the same number of samples, $\Nrisk\alpha=25$, lies in this tail for both targets, to maintain the relatively same error level for SAA of CVaR under different $\alpha$. The weight $\eta$ is heuristically tuned to ensure the empirical violation probability closely aligns with the target level $\alpha$.

In the following experiments, we compare the proposed LDT-CC method with the three baselines for the two risk measures in \Cref{sec:risk-measure}, using the following criteria.
We investigate the total iterations and time required for the optimization for each strategy, as well as the final objective value of the equilibrium cost \eqref{eq:J_ec_num}. Another important criterion we look at is the violation probability of the risk measure when evaluated using the optimized control $u^\ast$, i.e., $\PP(F(u^\ast, \theta)\geq z)$, which is empirically estimated using $10^6$ samples. With the same samples, we also report the density-weighted mean speed $\bar u := \mathbb{E}\big[\int_0^T\!\!\int_0^L \rho u\,dx\,dt \big/ \int_0^T\!\!\int_0^L \rho\,dx\,dt\big]$ and the throughput $\bar q_b := \mathbb{E}\big[\tfrac{1}{T}\int_0^T (\rho u)(x_b,t)\,dt\big]$ at location $x_b$. Here, $x_b$ is the first spatial grid point downstream of the source region $[0.4,0.6]$ in the normalized road coordinate.
Finally, all optimized controls are evaluated under a unified test driven by the exact MPP source $s(x,t;\theta^*)$. Case-specific settings are detailed in the following subsections.

\subsection{Risk control on maximum density}

For the risk measure of the maximum density, we set the safety threshold to $z = 0.65$ with an allowable violation probability of $\alpha = 0.5\%$ and $0.2\%$. The initial density is a Gaussian-like bottleneck: $\rho_0(x) = 0.06 + 0.54 \exp(-0.5 \sigma^{-2} (x - L/2)^2)$. We multiply the source \eqref{eq:kl} with a mask function $\ind{[0.4,0.6]\times[k, k+0.5]}(x,t), k\in \mathbb{Z}$, meaning that the spatio-temporal mask for the stochastic source $s(x,t;\theta)$ is constructed to be active within the spatial region $x \in [0.4, 0.6]$ and periodically active in time for the first $0.5\,\text{s}$ of every $1.0\,\text{s}$.

\begin{table}[tbp]
\vspace*{5pt}
\centering
\caption{Performance Comparison (Maximum Density)}
\label{tab:density_performance}
\scriptsize
\setlength{\tabcolsep}{1.5pt}
\renewcommand{\arraystretch}{1.1}
\begin{tabular}{l c c c c c c c}
\hline\hline
\textbf{Method} & \textbf{Target $\alpha$} & \textbf{Iter.} & \textbf{Time (s)} & \textbf{Viol.\ Prob.} & \textbf{Eq.\ Cost} & $\bar u$ & $\bar q_b$ \\
\hline
LDT-CC (Prop.) & 0.5\% & 174 & 6.25  & 0.53\% & 0.0009 & 0.533 & 0.253 \\
LDT-CC (Prop.) & 0.2\% & 217 & 7.43  & 0.25\% & 0.0009 & 0.533 & 0.254 \\
CVaR ($\eta=\!0.02252$) & 0.5\% & 33 & 123.85 & 0.58\% & 0.0012 & 0.533 & 0.255 \\
CVaR ($\eta=\!0.0223$) & 0.5\% & 19 & 73.78 & 52.96\% & 0.0008 & 0.529 & 0.248 \\
CVaR ($\eta=\!0.03162$) & 0.2\% & 41 & 483.67 & 0.22\% & 0.0016 & 0.531 & 0.256 \\
CVaR ($\eta=\!0.03$) & 0.2\% & 24 & 279.03 & 3.42\% & 0.0013 & 0.532 & 0.255 \\
Nonseparable & N/A & 31 & 0.79 & 99.19\% & 0.0043 & 0.545 & 0.261 \\
Unconstrained  & N/A & 25 & 0.89 & 100.00\% & 0.0002 & 0.525 & 0.243 \\
\hline\hline
\end{tabular}
\vspace{-5pt}
\end{table}

\begin{figure}[tbp]
    \centering
    \includegraphics[width=\linewidth]{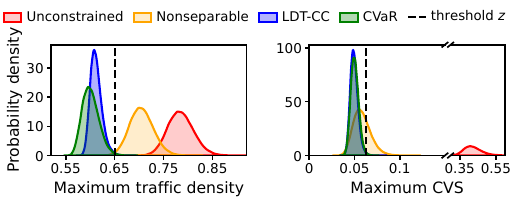}
    \caption{Probability density functions of the maximum density (left, $\alpha=0.5\%$) and of the maximum CVS (right) under the stochastic source term for the optimal controls obtained by the proposed method and the three baselines; the CVaR curves show the solutions with $\eta=0.02252$ (left) and $\eta=0.02433$ (right). The vertical dashed lines mark the thresholds $z=0.65$ and $z=0.063$. The broken x-axis in the right panel accommodates the different scales of risk of CVaR and LDT-CC (left part) and of the unconstrained strategy (right part).}
    \label{fig:kde_density}
\end{figure}

As shown in \Cref{tab:density_performance}, the unconstrained and nonseparable strategies lack direct risk measures, yielding violation probabilities near 100\% and completely failing to prevent rare density spikes. In contrast, the proposed LDT-CC achieves precise probability control, yielding violations of 0.53\% and 0.25\% for target levels $\alpha=0.5\%$ and 0.2\%, respectively, as \Cref{fig:kde_density} also shows. 
LDT-CC reaches the target $\alpha$ in one solve, while the CVaR approach needs to tune the weight $\eta$ for each target $\alpha$, by solving \eqref{eq:cvar_num} through a grid of $\eta$ and finding the optimal one whose violation probability is closest to $\alpha$.
With the selected optimal $\eta$, the CVaR method yields violation probabilities of $0.58\%$ and $0.22\%$ at equilibrium costs of $0.0012$ and $0.0016$, compared with a lower cost $0.0009$ achieved by LDT-CC. 
Note that the violation probability of the CVaR approach is extremely sensitive to the weight $\eta$. A change of one percent, from $\eta=0.02252$ to $\eta=0.0223$, raises it from $0.58\%$ to $52.96\%$.
Furthermore, \Cref{tab:density_performance} compares the computational cost. When the target probability $\alpha$ decreases from $0.5\%$ to $0.2\%$, the computation time of LDT-CC only slightly increases from $6.25\,\mathrm{s}$ to $7.43\,\mathrm{s}$, whereas the CVaR method (the optimal one over a grid of $\eta$) requires $123.85\,\mathrm{s}$ and $483.67\,\mathrm{s}$ for the reported solutions, since its number of risk samples grows in proportion to $1/\alpha$. 
LDT-CC takes about seven times the computation time of the unconstrained problem, and this is the cost of enforcing the probability target and is insensitive to the rareness $\alpha$. 
The mean speed and the throughput of LDT-CC are also reported in the last two columns and lie between those of the unconstrained and nonseparable strategies. 

\Cref{fig:kde_density} plots the probability density functions (PDFs) of the maximum traffic density. The unconstrained and nonseparable strategies strongly violate the safety threshold $z=0.65$ and span a large range, indicating high susceptibility to stochastic perturbations. While the CVaR method shifts the distribution near the target, it remains slightly wider than LDT-CC. In contrast, LDT-CC achieves the narrowest distribution and strictly bounds the most of the density below $z=0.65$.

\begin{figure}[tbp]
    \centering
    \includegraphics[width=1.0\linewidth]{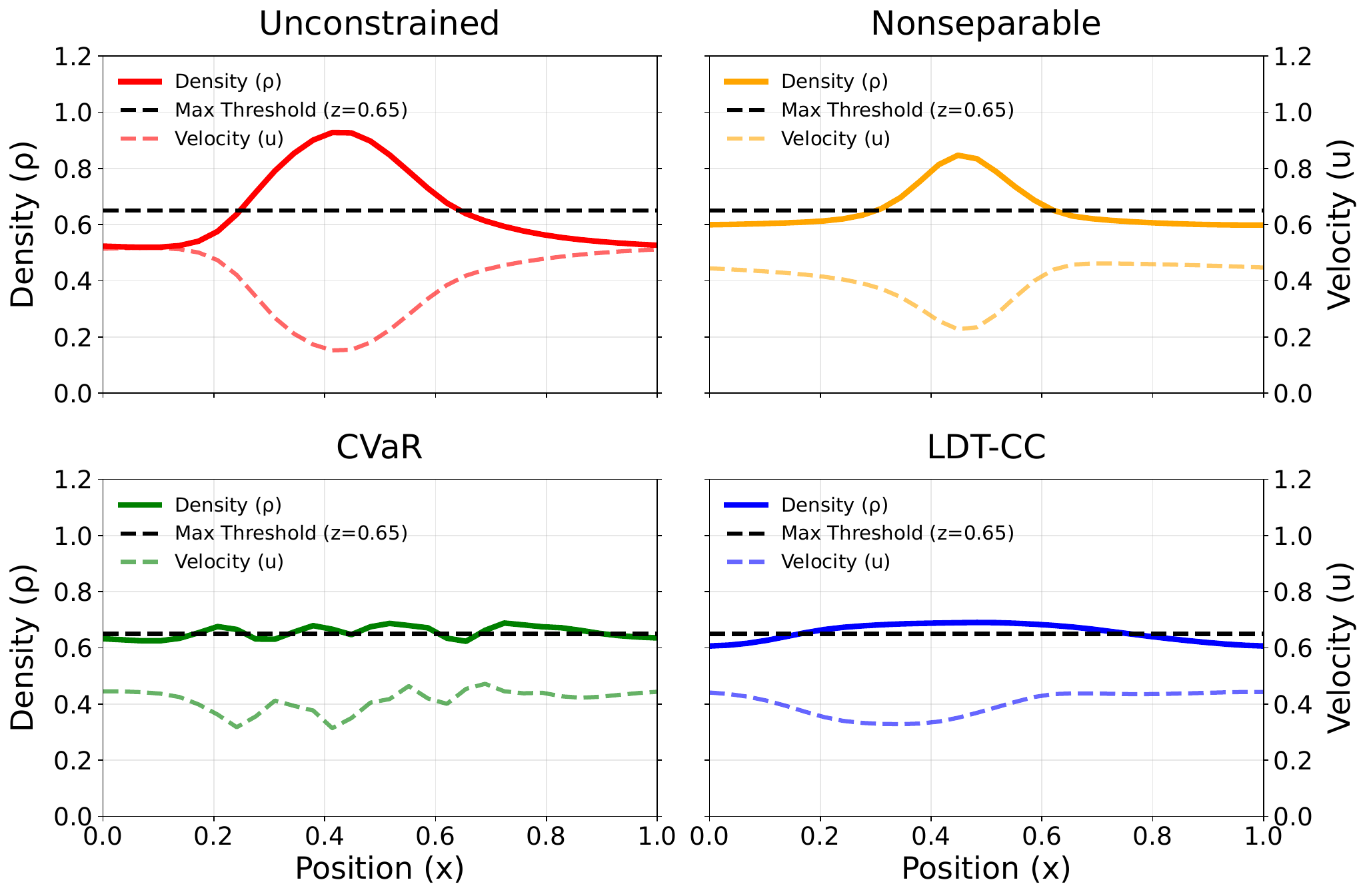}
    \caption{Snapshots at $t=4.5\,\mathrm{s}$ of the density and velocity for the maximum density case with $\alpha=0.2\%$ under the MPP source $s(x,t;\theta^\ast)$; the CVaR panel shows the solution with $\eta=0.03162$. The CVaR solution exhibits a sharp velocity field, whereas that of LDT-CC is smoother.
    }
    \label{fig:snap_shot_4.5s}

\end{figure}

\Cref{fig:snap_shot_4.5s} displays snapshots at $t=4.5\,\mathrm{s}$ of the traffic density and the optimal velocity $u^\ast$ under the MPP source $s(x,t;\theta^*)$ obtained from the LDT-CC method. The velocity field of the CVaR solution exhibits noticeable fluctuations in space, whereas that of LDT-CC is smooth. A smooth speed command is easier for the AVs to follow. Overall, the proposed optimal control satisfies the probability target while maintaining both a smooth velocity field and superior computational efficiency.

\subsection{Risk control on maximum CVS}

In the CVS-constrained case, we aim to limit local speed variation. We use a local spatial window of $\Delta=5$ grid cells. The threshold for the maximum CVS is $z=0.063$ with a target violation probability of $\alpha=0.5\%$. The system is initialized with a uniform density $\rho_0=0.2$, and the uncertain source is activated only during the initial interval via the mask  $\ind{[0.4,0.6]\times[0,2.0]}(x,t)$. 

As reported in \Cref{tab:cvs_performance}, the proposed LDT-CC achieves an empirical violation probability of $0.61\%$, remaining close to the $0.5\%$ target, and requires only $3.71\,\mathrm{s}$ of computation time. In contrast, due to the complexity of the CVS measure, the CVaR method requires significantly more computation time ($268.33\,\mathrm{s}$). The CVaR rows of \Cref{tab:cvs_performance} are solved from the same unconstrained solution over a grid of $\eta$. A change of one percent in the weight, from $\eta=0.02433$ to $\eta=0.0246$, raises the violation probability from $0.54\%$ to $99.99\%$. The weight controls the probability only indirectly, whereas the LDT-CC constraint \eqref{eq:nlp_prob} imposes it directly.

To further investigate these behaviors, \Cref{fig:kde_density} (right) plots the PDF of the maximum CVS. The unconstrained strategy yields a complete failure with a $100.00\%$ violation probability. This occurs because the CVS measure inherently uses local density $\rho(x,t)$ as a weight; without any safety constraints, the unconstrained strategy fails to prevent the simultaneous occurrence of severe speed heterogeneity and uneven density accumulation caused by the stochastic source. The nonseparable strategy reduces the violation probability to $28.86\%$ but still fails to satisfy the safety target. The reason here is different: the nonseparable policy naturally induces vehicles to decelerate in high-density areas and accelerate in low-density regions to track the equilibrium velocity. While this mechanism helps regulate traffic density, it actively generates local speed variation, which directly inflates the CVS measure. Overall, LDT-CC effectively enforces the complex speed variation limit with both high precision and superior computational efficiency, completely eliminating the need for heuristic parameter tuning.

\begin{table}[htbp]
\centering
\caption{Performance Comparison (Maximum CVS)}
\label{tab:cvs_performance}
\scriptsize
\setlength{\tabcolsep}{4pt}
\renewcommand{\arraystretch}{1.1}
\begin{tabular}{l c c c c c c}
\hline\hline
\textbf{Method} & \textbf{Iter.} & \textbf{Time (s)} & \textbf{Viol.\ Prob.} & \textbf{Eq.\ Cost} & $\bar u$ & $\bar q_b$ \\
\hline
LDT-CC (Prop.) & 95 & 3.71 & 0.61\% & 0.0012 & 0.513 & 0.259 \\
CVaR ($\eta=0.02433$) & 52 & 268.33 & 0.54\% & 0.0050 & 0.503 & 0.258 \\
CVaR ($\eta=0.0246$) & 28 & 165.14 & 99.99\% & 0.0012 & 0.508 & 0.255 \\
Nonseparable & 18 & 0.69 & 28.86\% & 0.0025 & 0.517 & 0.262 \\
Unconstrained  & 27 & 0.96 & 100.00\% & 0.0002 & 0.501 & 0.246 \\
\hline\hline
\end{tabular}
\end{table}

\section{Conclusions}

In this paper, we developed a sampling-free, rare chance-constrained optimal control framework (LDT-CC) for autonomous traffic management. By combining stochastic PDE models with an LDT-based approximation, we converted a challenging chance-constrained stochastic problem into a tractable deterministic NLP. This framework efficiently achieves precise near-target probability control for the maximum density and the maximum CVS risk measures. 
We tested its performance on a ring road with stochastic on-ramp inflow, demonstrating that our proposed LDT-CC requires no parameter tuning and maintains a low computational cost that is insensitive to event rarity.
Future work includes extensions to second-order macroscopic models, such as the Aw--Rascle--Zhang model, and to multiple sources of uncertainty, including boundary inflow fluctuations and random on-ramp demand.

\bibliographystyle{IEEEtran}
\bibliography{refs}

\end{document}